\documentclass[10pt]{article}

\usepackage{graphicx}
\usepackage{amsmath}
\usepackage{epsf}
\usepackage{graphpap}

\def\be{\begin{equation}}
\def\ee{\end{equation}}

\begin{document}

\begin{center}{\bf  \Large Closed-form solutions of Lucas-Uzawa model
with externalities via partial Hamiltonian approach }\\[2ex]
{R. Naz$^{a,*}$, Azam Chaudhry $^{b}$}\\[2ex]
{$^a$ Centre for Mathematics and Statistical Sciences,
  Lahore School of Economics, Lahore, 53200, Pakistan\\
  $^{b}$ Department of Economics,
Lahore School of Economics, Lahore, 53200, Pakistan\\
$^*$  Corresponding Author Email: drrehana@lahoreschool.edu.pk.\\

}
%
\end{center}
  \begin{abstract}
 In this paper, we establish multiple
closed-form solutions for all the variables in the Lucas-Uzawa model
with externalities for the case with no parameter restrictions as
well as for cases with specific parameter restrictions. These
multiple solutions are derived with the help of the results derived
in Naz et al (2016); Naz and Chaudhry (2017). This multiplicity of
solutions is new to the economic growth literature on Lucas-Uzawa
model with externalities. After finding solutions for the
Lucas-Uzawa model with externalities, we use these solutions to
derive the growth rates of all the variables in the system which
enables us to fully describe the dynamics of the model. The multiple
solutions can potentially explain why some countries economically
overtake other countries even though they start from the same
initial conditions.

\end{abstract}


\section{Introduction}

Over the last few decades, two strands of the literature looking at
economic growth have been brought together to explain the factors
that affect long run growth. So the neoclassical models of economic
growth, which were based on the assumption of competitive markets,
were adapted to incorporate the the- ories of human capital, which
looked at the private benefits of education. The first way in which
education was incorporated into economic growth theo- ries was by
incorporating human capital as a factor input which has a direct
impact on output (see Mankiw et al \cite{mank} and Romer
\cite{romer1}). The second way in which education was incorporated
into economic growth models was to include knowledge accumulation
into growth models through human capital accumu- lation or through
research and development activities (See Romer \cite{romer1,romer2},
Lucas \cite{lu} and Uzawa \cite{uz}).

What differentiated many of these models from traditional
neoclassical growth models was the way in which knowledge
accumulation was used to endogenize the factors that explain long
run economic growth.  These models extended the early growth models
which were characterized by exogenous productivity shocks and the
way endogeneity was introduced was through the presence of knowledge
externalities  (See Romer \cite{romer1}) in the economy-wide
production function which lead to increasing returns to scale.
Another way of introducing endogenous growth was developed by Lucas
\cite{lu} in the well-know Lucas-Uzawa model in which human capital
was incorporated into the production function in which the average
level of human capital in the economy affects individual firm level
output but is not taken into account in each firm's profit
maximizing decisions. Tamura \cite{tam} also developed a model of
economic growth which analyzed a human capital externality in the
production of human capital itself.

The basic Lucas-Uzawa model was solved using simplification methods:
Some authors have solved the basic Lucas-Uzawa model using dimension
reduction techniques (see Benhabib and Perli \cite{ben}, Caballe and
Santos \cite{cab}) or time elimination methods (see Mulligan and
Sala-i-Martin \cite{mul}), but these methods fail to provide
explicit solutions for the dynamics of the original models. Other
authors have used methods like parameter restrictions or
hypergeometric functions to find closed-form solutions of the
original Lucas-Uzawa model (see Xie \cite{xie}; Boucekkine and
Ruiz-Tamarit \cite{buk}). Recently Naz et al \cite{naz2016}
established closed-form solutions for the basic Lucas-Uzawa model by
using their newly developed partial Hamiltonian approach
\cite{naz,naz3} with a specific parameter restriction. The partial
Hamiltonian approach partial Hamiltonian approach \cite{naz,naz3} is
a significant contribution in Lie group theoretical methods. The
interested reader is referred to read some interesting papers
utilizing Lie group theoretical methods in solving real world
applications(e.g. \cite{imran1}-\cite{imran3}). Naz and Chaudhry
\cite{naz2017} then constructed completely new multiple closed-form
solutions for the basic Lucas-Uzawa model with no parameter
restrictions using partial Hamiltonian approach \cite{naz,naz3}.
Moreover,  Naz and Chaudhry \cite{naz2017} also compared the
closed-form solutions derived via partial Hamiltonian approach with
the solutions derived from the classical approach.

Ruiz-Tamarit \cite{ramon} have derived only one closed-form solution
of Lucas-Uzawa model with externalities under the assumption that
the capital share equals the inverse of the intertemporal elasticity
of substitution i.e. $\sigma=\beta$. Hiraguchi \cite{hir}
transformed the Lucas-Uzawa model with externalities into the basic
Lucas-Uzawa model and then derived a unique closed-form solution
directly from Boucekkine and Ruiz-Tamarit \cite{buk} for fairly
general values of the parameters of the model.

One important issue that arises is the role of growth models like
the Lucas-Uzawa model in explaining differences in cross country
growth rates especially between countries that start from the same
initial conditions. Though the Lucas-Uzawa model illustrates the key
role of human capital accumulation for endogenous growth, the basic
model may not explain how some countries overtake other countries
over time. This why the issue of multiplicity of equilibrium paths
is important in the context of growth models in general and the
Lucas-Uzawa model in particular.

This paper is distinct from the previous literature in that it
establishes multiple closed-form solutions for the Lucas-Uzawa model
with externalities.  This is the first time in the literature that
multiple closed-form solutions have been found for the Lucas-Uzawa
model with externalities and this multiplicity is critical in
understanding how certain countries may overtake other countries in
the growth process.  The closed-form solutions derived in this paper
are for the case with no parameter restrictions and also for the
cases with particular parameter restrictions.  We also derive the
growth rates of all the variables in the model and an interesting
feature of our result is that while one solution yields static
growth rates, the other solutions provide dynamic growth rates. In
the long run, all the growth rates approach the same steady state
growth rate.

The layout of the paper is as follows. In Section 2, we provide an
overview of the Lucas-Uzawa model with externalities based on the
transformation by Hiraguchi \cite{hir} . In Section 3, the
closed-form solutions of Lucas-Uzawa model with externalities are
derived for the no parameter restriction case.  The characteristics
of balanced growth path and growth rates of all the variables in the
model are analyzed in Section 4. In Section 5, we derive the
closed-form solutions of the model in the commonly discussed case
where $\sigma=\beta$. Finally, our conclusions are presented in
Section 6.

\section{Overview of Lucas-Uzawa model with externalities \cite{hir}}

In this section, we provide an overview of the transformed
Lucas-Uzwa model in which in Hiraguchi \cite{hir} reduced the
Lucas-Uzawa model with externalities to the basic model.

 \subsection{Model setup}
 The representative agent's utility function is
defined as
  \be Max_{c,k,h,u} \quad \int_0^{\infty}\frac{c^{1-\sigma}-1}{1-\sigma} e^{-\rho t} ,
  \; \sigma \not=1\label{(rn1)}\ee subject to constraints
on  physical capital and human capital: \begin{eqnarray} \dot k(t)
= \gamma k^\beta u^{1-\beta} h^{1-\beta+\theta}-\pi k-c, \; k_0=k(0) \nonumber\\
\dot h(t) =\delta(1-u)h,\; h_0=h(0) \label{(rn2)}\end{eqnarray}
where $1/{\sigma}$ is the constant elasticity of intertemporal
substitution, $\rho>0$ is the discount factor, $\beta$ is the
elasticity of output with respect to physical capital, $\gamma>0$ is
the technological levels in the good sector,  $\pi>0$ is the
depreciation rate for physical capital, $\delta>0$ is the
technological levels in the education sector, $k$ is the physical
capital, $h$ is the human capital, $c$ is per capita consumption and
$u$ is the fraction of labor allocated to the production of physical
capital. In the equation above, the parameter $\theta$ illustrates
the human capital externality (see Lucas (1988) and Tamura (1991)).
It is important to mention here that Hiraguchi \cite{hir} considered
a version of the Lucas-Uzwa model with externalities in which there
was no depreciation of physical capital i.e. $\pi=0$.

 The current value Hamiltonian function for this
problem is defined as \be
H(t,c,k,\lambda)=\frac{c^{1-\sigma}-1}{1-\sigma}+\lambda[\gamma
k^\beta u^{1-\beta} h^{1-\beta+\theta} -\pi k-c]+\mu\delta(1-u)h,
\label{(rn3)}\ee where $\lambda(t)$ and $\mu(t)$ are costate
variables.  The transversality conditions are \be \lim_{t\to\infty}
e^{-\rho t}\lambda(t)k(t)=0 ,\; \lim_{t\to\infty} e^{-\rho
t}\mu(t)h(t)=0
 \label{(trans)}. \ee

Pontrygin's maximum principle yields the following necessary first
order conditions for optimal
 control:
 \be
\lambda=c^{-\sigma}, \label{(rn4)} \ee \be
u^{\beta}=\frac{\gamma(1-\beta)k^{\beta}h^{-\beta+\theta}}{\delta}\frac{\lambda}{\mu},\label{(rn5)}
\ee
 \be
\dot k(t) =\gamma k^\beta u^{1-\beta} h^{1-\beta+\theta}-\pi k-c,
\label{(rn6)} \ee
 \be \dot h(t) =\delta(1-u)h, \label{(rn7)} \ee
 \be   \dot
\lambda= -\lambda \gamma \beta
u^{1-\beta}k^{\beta-1}h^{1-\beta+\theta}+\lambda (\rho+\pi),
\label{(rn8)}\ee
 \be   \dot
\mu= \mu(\rho-\delta)-\frac{\mu \theta \delta}{1-\beta}u
.\label{(rn9)}\ee

The growth rates for the variables $c$ and $u$ after simplification
take the following form: \be \frac{\dot c}{c}= \frac{ \beta
\gamma}{\sigma} u^{1-\beta}k^{\beta-1}h^{1-\beta+\theta}-\frac{
\rho+\pi}{\sigma}, \label{(rn10)} \ee \be \frac{\dot u}{u}=
\frac{(\delta+\pi)(1-\beta)+\delta \theta}{\beta}-\frac{ c}{k}+
\bigg (\frac{1-\beta+\theta}{1-\beta} \bigg) \delta u.
\label{(rn11)} \ee

\subsection{Transformed first-order conditions \cite{hir}}
Introducing the transformations, put forward by Hiraguchi
\cite{hir}, \be \phi=\frac{1-\beta+\theta}{1-\beta},\;
h^*=h^{\phi},\; \delta^*=\delta \phi,\; \mu^*=\mu \phi^{-1}
h^{1-\phi},\label{(trans)}\ee the system of equations {(\ref
{(rn4)})}-{(\ref {(rn11)})} takes the following form:
 \be
\lambda=c^{-\sigma}, \label{(bb4)} \ee \be \bigg(\frac{h^*
u}{k}\bigg)^{\beta}=\frac{\gamma(1-\beta)}{\delta^*}\frac{\lambda}{\mu^*},\label{(bb5)}
\ee
 \be
\dot k=\gamma k^\beta(u h^*)^{1-\beta} -\pi k-c,\label{(bb6)} \ee
 \be \dot h^*(t) =\delta^*(1-u)h^* ,\label{(bb7)} \ee
 \be   \dot
\lambda= -\lambda \beta \gamma \bigg(\frac{h^*
u}{k}\bigg)^{1-\beta}+\lambda (\rho+\pi), \label{(bb8)}\ee
 \be   \dot
\mu^*= \mu^*(\rho-\delta^*),\label{(bb9)}\ee \be \frac{\dot c}{c}=
\frac{ \beta \gamma}{\sigma} \bigg(\frac{h^*
u}{k}\bigg)^{1-\beta}-\frac{ \rho+\pi}{\sigma}, \label{(bb10)} \ee
\be \frac{\dot u}{u}= \frac{(\delta^*+\pi)(1-\beta)}{\beta}-\frac{
c}{k}+\delta^* u. \label{(bb11)} \ee
 The transversality conditions transform to \be \lim_{t\to\infty}
e^{-\rho t}\lambda(t)k(t)=0 ,\; \lim_{t\to\infty} e^{-\rho
t}\mu^*(t)h^*(t)=0
 \label{(rnTC)}. \ee

 \subsection{Transformed model \cite{hir}}
Based on the transformation discussed above, Hiraguchi \cite{hir}
considered the following new growth problem with the variable $h^*$:

 \be Max_{c,k,h^*,u} \quad \int_0^{\infty}\frac{c^{1-\sigma}-1}{1-\sigma} e^{-\rho t} ,
  \; \sigma \not=1\label{(rn1new)}\ee subject to constraints
on  physical capital and human capital: \begin{eqnarray} \dot k(t)
= \gamma k^\beta u^{1-\beta} h^{*1-\beta+\theta}-\pi k-c, \; k_0=k(0) \nonumber\\
\dot h(t) =\delta^*h^*(1-u)h,\; h^*_0=h^*(0),
\label{(rn2new)}\end{eqnarray} where
$\phi=\frac{1-\beta+\theta}{1-\beta},\; h^*=h^{\phi},\;
\delta^*=\delta \phi.$

 The current value Hamiltonian for this problem is
 \begin{eqnarray}
H^*(t,c,u,k,h^*,\lambda,\mu^*)=\frac{c^{1-\sigma}-1}{1-\sigma}+\lambda[\gamma
k^\beta(u h^*)^{1-\beta} -\pi k-c]\nonumber\\+\mu^*\delta^*(1-u)h^*,
\label{(bb3new)}\end{eqnarray} where $\mu^*=\mu \phi^{-1}
h^{1-\phi}$.  It is straight forward to show that the first order
conditions and transversality conditions of this new problem are the
same as given in {(\ref {(bb4)})}-{(\ref {(rnTC)})}.

  Let $\bar{c},\bar{u},\bar{k},\bar{h}^*$ be the equilibrium
values for consumption $c$, the fraction of labor allocated to the
production of physical capital $u$, physical capital $k$ and human
capital $h^*$. Let $g_c, g_k, g_{h^*}$ and $g_u$ be the growth rates
 of the per capita consumption $c$,
physical capital $k$, human capital $h^*$, the fraction of labor
allocated to the production of physical capital $u$. In long run,
for $\rho<\delta^* <\rho+\delta^* \sigma$, the balanced growth path
(BGP) satisfies the following conditions stated (see e.g. Naz et al (2016); Naz and Chaudhry (2017)): \\

\begin{eqnarray}
g_c=g_k=g_{h^*}=\frac{\delta^*-\rho}{\sigma},\;g_u=0\nonumber\\
\bar{u}=\frac{\rho+\delta^*(\sigma-1)}{\delta \sigma},\; \bar{u}\in
[0,1],\nonumber\\
\frac{\bar{c}}{\bar{k}}=\frac{\delta^*+\pi(1-\beta)}{\beta}-\frac{\delta^*-\rho}{\sigma}=\xi>0,\label{(bgp)}\\
\frac{\bar{k}}{\bar{h^*}}=\frac{\rho+\delta^*(\sigma-1)}{\delta^*
\sigma}\bigg(\frac{\beta \gamma}{\delta^*+\pi}
\bigg)^{\frac{1}{1-\beta}}.\nonumber
\end{eqnarray}

The characteristics of the BGP for the Lucas-Uzawa model with
externalities are given in following proposition:\\

{\bf Proposition 1:} Let $g_c, g_k, g_h$ and $g_u$ be the growth
rates
 of the per capita consumption $c$,
physical capital $k$, human capital $h$, the fraction of labor
allocated to the production of physical capital $u$. Let
$\rho(1-\beta)<\delta(1-\beta+\theta)<
\rho(1-\beta)+\delta\sigma(1-\beta+\theta)$, then the system reaches the BGP and the following statements are valid:\\

 i. $ g_c=g_k=\frac{(\delta
-\rho)(1-\beta)+\delta \theta}{\sigma(1-\beta)}$, $g_h=\frac{(\delta
-\rho)(1-\beta)+\delta
\theta}{\sigma(1-\beta+\theta)}$ and $g_u=0$,\\
ii. $\bar{u}\in [0,1]$ and $\bar{u}=\frac{(\rho-\delta+\sigma
\delta)(1-\beta)+\delta \theta (\sigma-1)}{\delta \sigma
(1-\beta+\theta)}$,\\
iii.
$\frac{\bar{c}}{\bar{k}}=\frac{\delta(1-\beta+\theta)+\pi(1-\beta)^2}{\beta(1-\beta)}-\frac{\delta(1-\beta+\theta)-\rho(1-\beta)}{\sigma
(1-\beta)}=\xi$\\
iv. $\frac{\bar{k}}{\bar{h^{\phi}}}=
\frac{\bar{u}}{\bar{z}}=\frac{(\rho-\delta+\delta
\sigma)(1-\beta)+\delta \theta(\sigma-1)}{\delta \sigma
(1-\beta+\theta)}\bigg(\frac{\beta \gamma
(1-\beta)}{\delta(1-\beta+\theta) +\pi(1-\beta)}
\bigg)^{\frac{1}{\beta-1}}$.\\

 {\bf Proof}: The results directly follow by using  $\phi=\frac{1-\beta+\theta}{1-\beta},\; h^*=h^{\phi},\;
\delta^*=\delta \phi$ in {(\ref {(bgp)})}. We need to prove $\xi>0$.
 \begin{eqnarray} \xi=
\frac{\delta(1-\beta+\theta)+\pi(1-\beta)^2}{\beta(1-\beta)}-\frac{\delta(1-\beta+\theta)-\rho(1-\beta)}{\sigma
(1-\beta)}\nonumber\\
=\frac{\sigma \delta(1-\beta+\theta)+\sigma\pi(1-\beta)^2-\beta
\delta(1-\beta+\theta)+\beta \rho(1-\beta)}{\sigma
\beta(1-\beta)}\nonumber\\
>\frac{[\delta(1-\beta+\theta)-\rho(1-\beta)])[\delta(1-\beta+\theta)+\pi(1-\beta)])(1-\beta)}{\delta \sigma
\beta(1-\beta+\theta)}>0
\end{eqnarray}
and this completes the proof.\\

In the next section, we establish the closed-form solutions of
Lucas-Uzawa model with externalities by utilizing the results of
basic Lucas-Uzawa model presented in \cite{naz2016,naz2017}. Then we
check all the derived closed-form solutions to ensure that they
satify all the properties of the BGP.

\section{Closed-form solutions for Lucas-Uzawa-model with
externalities}  In this section, we derive multiple closed-form
solutions for the Lucas-Uzawa model with externalities by utilizing
the results of the basic Lucas-Uzawa model presented in
\cite{naz2016,naz2017}. If we let $[k,h^*,c,u,\lambda,\mu^*]$ denote
the closed-form solution to the transformed problem, a closed-form
solution for the original model also exists and is given by
$[k,h,c,u,\lambda,\mu]=[k,(h^*)^{\frac{1}{\phi}},c,u,\lambda,\phi
\mu^* (h^*)^{\phi-1}]$  (see e.g. \cite{hir}).
\subsection{Closed-form solutions using fairly general parameters values}
Naz and Chaudhry \cite{naz2017} established three sets of
closed-from solutions for the basic Lucas-Uzawa model with the help
of their newly developed partial Hamiltonian approach for fairly
general values of the parameters. We derive a solution for the
dynamical system of ODEs {(\ref {(bb4)})}-{(\ref {(bb11)})} directly
from \cite{naz2017} in terms of a variable
$z(t)=\frac{u(t)h^*(t)}{k(t)}$.

The first set of solutions for all the economic variables of the
dynamical system of ODEs {(\ref {(bb4)})}-{(\ref {(bb11)})},
satisfying the transversality condition {(\ref {(rnTC)})}, is (see
\cite{naz2017}):
\begin{eqnarray} c(t)=c_0e^{-\frac{(\rho-\delta^*
)}{\sigma}t}, \nonumber\\
 k(t)=k_0e^{-\frac{(\rho-\delta^* )}{\sigma}t},\nonumber\\
  u(t)=\frac{\rho-\delta^*+\delta^*
\sigma}{\delta^* \sigma}=\bar{u},\label{(sol1)} \\
 h^*(t)=h^*_0 e^{-\frac{(\rho-\delta^* )}{\sigma
}t}, \nonumber\\
\lambda(t)=c_0^{-\sigma}e^{(\rho-\delta^*)t},\nonumber\\
 \mu^*= c_1  e^{(\rho-\delta^*)
 t},\nonumber\\
 \bar{z}=\bigg(\frac {\delta^* +\pi}{\beta \gamma}
  \bigg)^{\frac{1}{1-\beta}} \nonumber,\end{eqnarray}
  provided $ \rho<\delta^* <\rho+\delta^* \sigma ,
c_0=\bigg(\frac{(1-\beta)\gamma }{c_1 \delta^* \bar{z}^{ \beta}
}\bigg)^{\frac{1}{\sigma}},
 \frac{c_0}{k_0}=\frac{\delta^*+\pi(1-\beta)}{\beta}-\frac{\delta^*-\rho}{\sigma},h^*_0=\frac{\bar{z} k_0}{\bar{u}}$.

 The second set of solutions for all the economic variables of the
dynamical system of ODEs {(\ref {(bb4)})}-{(\ref {(bb11)})},
satisfying the transversality condition {(\ref {(rnTC)})}, is (see
\cite{naz2017}):
\begin{eqnarray} c(t)=c_0
z_0^{\frac{\beta}{\sigma}}e^{-\frac{(\rho-\delta^*
)}{\sigma}t}z^{-\frac{\beta}{\sigma}}, \nonumber\\
k(t)=\bigg(\frac{k_0}{c_0z_0^{\frac{\beta-\sigma}{\sigma}}}-F(t)
\bigg) c_0z_0^{\frac{\beta}{\sigma}}z(t)^{-1}e^{\frac{(\delta^*
+\pi-\pi \beta)}{\beta}t}
,\nonumber\\
 h^*(t)=\frac{h^*_0}{z_0[\sigma c_0 z_0^{\beta-1}
-(\rho+\pi-\pi \sigma)k_0 z_0^{\beta-1}+\beta \gamma(1-\sigma)k_0]}[
\sigma c_0 z_0^{\frac{\beta}{\sigma}}e^{-\frac{(\rho-\delta^*
)}{\sigma}t}z^{-\frac{\beta}{\sigma}+\beta} \nonumber\\
+(\beta \gamma (1-\sigma) -(\rho+\pi-\pi \sigma)z^{\beta-1})
(\frac{k_0}{c_0z_0^{\frac{\beta-\sigma}{\sigma}}}-F(t) )
c_0z_0^{\frac{\beta}{\sigma}}e^{\frac{(\delta^*+\pi
-\pi \beta)}{\beta}t}], \nonumber\\
 u(t)=\frac{u_0}{k_0}[\sigma c_0 z_0^{\beta-1}-(\rho+\pi-\pi \sigma)k_0 z_0^{\beta-1}+\beta \gamma(1-\sigma)k_0]\nonumber\\
 \times \frac{(\frac{k_0}{c_0
z_0^{\frac{\beta-\sigma}{\sigma}}}-F(t))}{[\beta
\gamma(1-\sigma)-(\rho+\pi-\pi \sigma)
 z^{\beta-1}](\frac{k_0}{c_0 z_0^{\frac{\beta-\sigma}{\sigma}}}-F(t))+\sigma z^{\beta-\frac{\beta}{\sigma}}e^{-(\frac{\delta^*+\pi-\pi
\beta}{\beta}-\frac{\delta^*-\rho}{\sigma} )t}},\nonumber\\
 \lambda(t)=c_0^{-\sigma}z_0^{-\beta}e^{(\rho-\delta^*)t}z^{\beta},\nonumber\\
 \mu^*(t)= c_1 e^{(\rho-\delta^* )t}
\nonumber,\end{eqnarray} where \begin{eqnarray} F(t)=\int_0^t
z(t)^{\frac{\sigma-\beta}{\sigma}} e^{-(\frac{\delta^*+\pi-\pi
\beta}{\beta}-\frac{\delta^*-\rho}{\sigma} )t} dt,\nonumber\\ z(t)=
\frac{\bar{z}z_0}{[(\bar{z}^{1-\beta}-z_0^{1-\beta})e^{-\frac{(1-\beta)(\delta^*
+\pi)}{\beta}t}+z_0^{1-\beta}]^{\frac{1}{1-\beta}} },\label{(sol2)}\\
\lim_{t\to\infty}
F(t)=\frac{k_0}{c_0z_0^{\frac{\beta-\sigma}{\sigma}}},\nonumber\\
 \rho<\delta^* <\rho+\delta^* \sigma, \xi=\frac{\delta^*+\pi-\pi
\beta}{\beta}-\frac{\delta^*-\rho}{\sigma}
>0,\nonumber\\
c_0 z_0^{\frac{\beta}{\sigma}}=\bigg(\frac{c_1
\delta^*}{(1-\beta)\gamma}\bigg)^{-\frac{1}{\sigma}},\nonumber\\
\frac{\gamma(1-\beta)(\rho-\delta^*+\delta^* \sigma)}{\delta^*}\nonumber\\=\frac{u_0}{k_0}[\sigma c_0 z_0^{\beta-1}-(\rho+\pi-\pi \sigma)k_0 z_0^{\beta-1}+\beta \gamma(1-\sigma)k_0],\nonumber\\
 \bar{z}=\bigg(\frac{\beta \gamma}{\delta^* +\pi} \bigg
)^{\frac{1}{\beta-1}}\nonumber.\end{eqnarray}

 The third set of solutions for all the economic variables of the
dynamical system of ODEs {(\ref {(bb4)})}-{(\ref {(bb11)})},
satisfying the transversality condition {(\ref {(rnTC)})}, is (see
\cite{naz2017}):
\begin{eqnarray} c(t)=c_0
z_0^{\frac{\beta}{\sigma}}e^{-\frac{(\rho-\delta^*
)}{\sigma}t}z^{-\frac{\beta}{\sigma}}, \nonumber\\
k(t)=\bigg(\frac{k_0}{c_0z_0^{\frac{\beta-\sigma}{\sigma}}}-F(t)
\bigg) c_0z_0^{\frac{\beta}{\sigma}}z(t)^{-1}e^{\frac{(\pi+\delta^*
-\pi \beta)}{\beta}t}
,\nonumber\\
h^*(t)=\bigg[\bigg(\frac{(\delta^*+\pi)(1-\beta)}{\beta} \frac{k_0
}{c_0 z_0^{\frac{\beta-\sigma}{\sigma}} } +\frac{\delta^* u_0 k_0
}{c_0 z_0^{\frac{\beta-\sigma}{\sigma}} }-\delta^* u_0
G(t)\bigg)e^{ -\frac{(\delta^*+\pi)(1-\beta)}{\beta}t}\nonumber\\
-\delta^* u_0 (\frac{k_0 }{c_0 z_0^{\frac{\beta-\sigma}{\sigma}}}-
F(t))\bigg]\times
\frac{c_0z_0^{\frac{\beta}{\sigma}}}{\frac{(\delta^*+\pi)(1-\beta)}{\beta}
u_0}e^{\frac{(\pi+\delta^* -\pi \beta)}{\beta}t},
\nonumber \\
u(t)=\frac{\frac{(\delta^*+\pi)(1-\beta)}{\beta} u_0[\frac{k_0 }{c_0
z_0^{\frac{\beta-\sigma}{\sigma}}}-
F(t)]}{[(\frac{(\delta^*+\pi)(1-\beta)}{\beta}  +\delta^*
u_0)\frac{k_0 }{c_0 z_0^{\frac{\beta-\sigma}{\sigma}} }-\delta^* u_0
G(t)]e^{ -\frac{(\delta^*+\pi)(1-\beta)}{\beta}t}-\delta^* u_0
[\frac{k_0 }{c_0
z_0^{\frac{\beta-\sigma}{\sigma}}}- F(t)]}\nonumber\\
 \lambda(t)=c_0^{-\sigma}z_0^{-\beta}e^{(\rho-\delta^*)t}z^{\beta},\nonumber\\
 \mu^*(t)= c_1 e^{(\rho-\delta^* )t}
\nonumber,\end{eqnarray}  where

\begin{eqnarray} \rho<\delta^* <\rho+\delta^* \sigma, \frac{\delta^*+\pi-\pi
\beta}{\beta}-\frac{\delta^*-\rho}{\sigma}
>0,\nonumber\\ F(t)=\int_0^t
z(t)^{\frac{\sigma-\beta}{\sigma}} e^{-(\frac{\delta^*+\pi-\pi
\beta}{\beta}-\frac{\delta^*-\rho}{\sigma} )t} dt, \nonumber\\
G(t)=\int_0^t z(t)^{\frac{\sigma-\beta}{\sigma}} e^{-\frac{\delta^*
\sigma-\delta^*+\rho}{\sigma} t} dt,
\label{(sol3)} \\
z(t)=
\frac{\bar{z}z_0}{[(\bar{z}^{1-\beta}-z_0^{1-\beta})e^{-\frac{(1-\beta)(\delta^*
+\pi)}{\beta}t}+z_0^{1-\beta}]^{\frac{1}{1-\beta}} },\nonumber\\
 c_0
z_0^{\frac{\beta}{\sigma}}=\bigg(\frac{c_1
\delta^*}{(1-\beta)\gamma}\bigg)^{-\frac{1}{\sigma}}, \nonumber\\
\lim_{t\to\infty}
F(t)=\frac{k_0}{c_0z_0^{\frac{\beta-\sigma}{\sigma}}},\nonumber\\
\lim_{t\to\infty} \bigg[(\frac{(\delta^*+\pi)(1-\beta)}{\beta}
+\delta^* u_0)\frac{k_0 }{c_0 z_0^{\frac{\beta-\sigma}{\sigma}}
}-\delta^* u_0 G(t) \bigg]=0,\nonumber\\ \lim_{t\to\infty} G(t)=
\frac{(\frac{(\delta^*+\pi)(1-\beta)}{\beta}  +\delta^*
u_0)}{\delta^* u_0}\lim_{t\to\infty} F(t),
\nonumber\\\bar{z}=\bigg(\frac{\beta \gamma}{\delta^* +\pi} \bigg
)^{\frac{1}{\beta-1}}.\nonumber\end{eqnarray}  The closed-form
solutions {(\ref {(sol1)})},{(\ref {(sol2)})} and {(\ref {(sol3)})}
can be expressed in terms of the original variables of the model by
using the transformations $ \phi=\frac{1-\beta+\theta}{1-\beta},\;
h^*=h^{\phi},\; \delta^*=\delta \phi,\; \mu^*=\mu \phi^{-1}
h^{1-\phi}.$ This completes the closed-form solutions of Lucas-Uzawa
model with externalities with no parameter restrictions.

\subsection{Comparison of the multiple closed-form solutions} Since
this is the first time in the literature that multiple closed-form
solutions have been obtained for the Lucas-Uzawa model with
externalities, it is useful to compare these results.

The first thing to note is that the values of consumption $c$,
physical capital stock $k$, and the costate variables $\lambda$ and
$\mu$ in solution {(\ref {(sol1)})} are different from the values of
these variables in solutions {(\ref {(sol2)})} and {(\ref
{(sol3)})}. Next, a comparison of the closed-form solutions {(\ref
{(sol2)})} and {(\ref {(sol3)})}  shows that the expressions for
consumption $c$, physical capital stock $k$, and the costate
variables $\lambda$ and $\mu$ are the same in both solutions. On the
other hand, the expressions for the fraction of labor devoted to
physical capital, $u$, and the level of human capital, $h^*$, are
different.

These results are important in the context of understanding
differences between economies transitioning towards their long run
equilibria. Our results show that countries that start with the same
initial conditions can have significant differences as they progress
towards their long run equilibria. For example, our results show
that countries can have the same levels of consumption and capital
stock, but may differ significantly in their levels of human
capital. Or countries may have the same levels of consumption but
may differ significantly in their levels of physical and human
capital. These differences are critical in understanding differences
between countries during the development process.

In terms of the previous literature, Hiraguchi \cite{hir} derived
only one solution for the model without parameter restrictions which
was similar to {(\ref {(sol3)})}, and in this solution $F(t)$ and
$G(t)$ were computed in terms of the hypergeometric functions.
Hiraguchi \cite{hir} also claimed that there existed no other
solutions for the Lucas-Uzawa model with externalities.  Contrary to
this claim and for the first time in the literature, we have
established multiple closed-form solutions for the Lucas-Uzwa model
with externalities in the case where there are no parameter
restrictions.

\section{Characteristics of balanced growth path}
In this section, we discuss the growth rates of the key variables in
our model for each of the solutions we have obtained. We then
compare the long run equilibrium values of these growth rates.

\subsection{Growth rates for each of the solutions}

First of all, we analyze the simple solution {(\ref {(sol1)})}. The
growth rates for per capita consumption $c$, physical capital $k$,
human capital $h$, the fraction of labor allocated to the production
of physical capital $u$, costate variables $\mu$ and $\lambda$ for
solution
 {(\ref
{(sol1)})} after using
 $
\phi=\frac{1-\beta+\theta}{1-\beta},\; h^*=h^{\phi},\;
\delta^*=\delta \phi,\; \mu^*=\mu \phi^{-1} h^{1-\phi},$  (and after
some simplifications) take the following forms:
\begin{eqnarray}
 \frac{\dot c}{c}=\frac{(\delta -\rho)(1-\beta)+\delta \theta}{\sigma(1-\beta)}, \nonumber \\
\frac{\dot k}{k}=\frac{(\delta -\rho)(1-\beta)+\delta
\theta}{\sigma(1-\beta)},\nonumber \\
 \frac{\dot
{h}}{h}=\frac{(\delta -\rho)(1-\beta)+\delta
\theta}{\sigma(1-\beta+\theta)},\nonumber \\
 \frac{\dot
u}{u}=0,\label{(gr100)}\\
 \frac{\dot
{\lambda}}{\lambda}=\frac{(\rho-\delta)(1-\beta)-\delta
\theta}{1-\beta},\nonumber \\
 \frac{\dot
{\mu}}{\mu}=\frac{((\rho-\delta)(1-\beta)-\delta \theta
)(\sigma(1-\beta+\theta)-\theta)}{\sigma(1-\beta)(1-\beta+\theta)}.
\nonumber
\end{eqnarray}

Solution {(\ref {(sol1)})} yields constant growth rates {(\ref
{(gr100)})} for all the variables of the model. The growth rates of
the per capita consumption $c$, physical capital $k$ and human
capital $h$ are positive provided $(\delta -\rho)(1-\beta)+\delta
\theta>0$.\\

The growth rates of all the variables for closed-form solution
{(\ref {(sol2)})} after simplifications are {\small
\begin{eqnarray}
 \frac{\dot c}{c}=\frac{(\delta-\rho)(1-\beta)+\delta \theta}{\sigma(1-\beta)}-\frac{\beta}{\sigma} \frac{\dot z}{z}, \nonumber \\
\frac{\dot
k}{k}=\frac{\pi(1-\beta)^2+\delta(1-\beta+\theta)}{\beta(1-\beta)}-\frac{z(t)^{\frac{\sigma-\beta}{\sigma}}
e^{-\xi t}}{\frac{k_0
}{c_0 z_0^{\frac{\sigma-\beta}{\beta}}}-F(t)}-\frac{\dot z}{z},\nonumber \\
 \frac{\dot
{h}}{h}=\frac{1-\beta}{1-\beta+\theta}\bigg[\frac{\sigma z(t)^{\beta-1}}{\sigma z(t)^{\beta-1}+[\beta \gamma(1-\sigma)-(\rho+\pi-\pi \sigma)z(t)^{\beta-1}]\frac{k(t)}{c(t)}}\frac{\dot c}{c}\nonumber \\
+\frac{\sigma \beta z(t)^{\beta-1}
\frac{c(t)}{k(t)}-(\rho+\pi-\pi\sigma)(\beta-1)z(t)^{\beta-1}+[\beta
\gamma(1-\sigma)-(\rho+\pi-\pi \sigma)z(t)^{\beta-1}]}{\sigma
z(t)^{\beta-1}\frac{c(t)}{k(t)} +\beta
\gamma(1-\sigma)-(\rho+\pi-\pi \sigma)z(t)^{\beta-1}}\frac{\dot
z}{z}\nonumber\\
+ \frac{\beta \gamma(1-\sigma)-(\rho+\pi-\pi
\sigma)z(t)^{\beta-1}}{\sigma z(t)^{\beta-1}\frac{c(t)}{k(t)} +\beta
\gamma(1-\sigma)-(\rho+\pi-\pi \sigma)z(t)^{\beta-1}}\frac{\dot
k}{k}\bigg]\nonumber\\
 \frac{\dot
u}{u}=\frac{\dot
k}{k}-\bigg(\frac{1-\beta+\theta}{1-\beta}\bigg)\frac{\dot
{h}}{h}+\frac{\dot z}{z},\label{(gr7)}\\
 \frac{\dot
{\lambda}}{\lambda}=\frac{(\rho-\delta)(1-\beta)-\delta \theta}{(1-\beta)}+\beta \frac{\dot z}{z},\nonumber \\
 \frac{\dot{\mu}}{\mu}=\frac{((\rho-\delta)(1-\beta)-\delta \theta
)(\sigma(1-\beta+\theta)-\theta)}{\sigma(1-\beta)(1-\beta+\theta)},
\nonumber\\
\xi=\frac{\delta^*+\pi-\pi
\beta}{\beta}-\frac{\delta^*-\rho}{\sigma} \nonumber
\end{eqnarray}}

where \be \frac{\dot
z}{z}=\frac{(z_0^{1-\beta}-\bar{z}^{1-\beta})\gamma
\bar{z}^{1-\beta} { {\rm e}^{{-\frac { \left( 1-\beta \right) \left(
\pi +\delta \right) +\delta \theta }{\beta}t}}} }{
\bar{z}^{\beta-1}+(z_0^{1-\beta}-\bar{z}^{1-\beta}){{\rm e}^{{-\frac
{ \left( 1-\beta \right)
 \left( \pi +\delta \phi \right)+\delta \theta }{\beta}}t}}
}.\label{(grz)} \ee

It can be shown that for solution {(\ref {(sol3)})} the dynamic
growth rates for consumption $c$, physical capital stock $k$, the
fraction of labor allocated to the production of physical capital
$u$ and the costate variables $\lambda$ and $\mu$ will be also equal
to the expressions for these variables given in {(\ref {(gr7)})}.
The growth rate for human capital human capital $h$ for solution
{(\ref {(sol3)})} is complex and is omitted.

\subsection{Comparison of the growth rates}

First, it is important to point out that we obtain a static growth
rate {(\ref {(gr100)})} for one solution  and dynamic growth rates
 {(\ref {(gr7)})}  for the other solutions. Next,
it is important to see what happens to these growth rates in the
long run.

In order to determine the long run growth rates, once needs to start
by looking at the long run value of $\frac{\dot z}{z}$.  It is clear
from equation {(\ref {(grz)})} that $\frac{\dot z}{z}$ approaches
zero as $t\mapsto\infty$ which means that the rate of growth of $z$
decreases asymptotically as we approach the steady state.

This means that in the long run: (i) the growth rates of consumption
$c$ and the physical capital $k$ decrease over time and approach
$\frac{(\delta -\rho)(1-\beta)+\delta \theta}{\sigma(1-\beta)}$ as
$t\mapsto\infty$; (ii) the growth rate of human capital decreases
over time and approaches $\frac{(\delta -\rho)(1-\beta)+\delta
\theta}{\sigma(1-\beta+\theta)}$  as $t\mapsto\infty$; (iii) the
growth rate of the fraction of labor allocated to the production of
physical capital $u$ approaches zero as $t\mapsto\infty$; (iv) The
growth rate of costate variable $\lambda$ equals
$\frac{(\rho-\delta)(1-\beta)-\delta \theta}{1-\beta}$ as
$t\mapsto\infty$; and (v) The growth rate of costate variable $\mu$
equals $ \frac{((\rho-\delta)(1-\beta)-\delta \theta
)(\sigma(1-\beta+\theta)-\theta)}{\sigma(1-\beta)(1-\beta+\theta)}$
as $t\mapsto\infty$.

Our results imply that even though the growth rates of the key
variables in the model differ between solutions, in the long run,
the static growth rates of all the variables in {(\ref {(gr100)})}
and the dynamic growth rates of all the variables in {(\ref
{(gr7)})} reach the same value. Or, in other words, even though the
key variables may differ in terms of growth between countries in the
short run, in the long run all the economies reach the same steady
state growth rates.

In order to confirm that our solutions satisfy the conditions for
the BGP (given in Proposition 1), it is also straight forward to
show using l'H$\hat{o}$pital rule that for the closed-form solutions
{(\ref {(sol1)})}, {(\ref {(sol2)})} and {(\ref {(sol3)})}, that \be
\lim_{t\to\infty} u(t)=\bar{u}.\ee

Also, for the closed-form solutions {(\ref {(sol1)})}, {(\ref
{(sol2)})} and {(\ref {(sol3)})} the ratios $\frac{c(t)}{k(t)}$ and
$\frac{k(t)}{h^{\phi}(t)}$ equal the values given in Proposition 1.
For the sake of brevity, we are omitting these calculations.

 \section{Lucas-Uzawa model with externalities for the $\sigma=\beta$ case}
 In this Section, we consider the special case when $\sigma=\beta$ which is a simplifying assumption that is  commonly used in the economic growth
 literature.  We start by deriving the closed-form solutions for this special case and
 then discuss the balanced growth path associated with this
 solution.

\subsection{Closed-from solution for the $\sigma=\beta$ case}
 We obtain two closed-form solutions for the Lucas-Uzawa model with externalities under the restriction
$\sigma=\beta$. The first closed-form solution for all variables is
obtained by taking $\sigma=\beta$ in {(\ref {(sol1)})} and also
using the transformations
 $
\phi=\frac{1-\beta+\theta}{1-\beta},\; h^*=h^{\phi},\;
\delta^*=\delta \phi,\; \mu^*=\mu \phi^{-1} h^{1-\phi},$ and is
given by

  \begin{eqnarray}  c(t)=c_0e^{-\frac{(\rho-\delta)(1-\beta)-\delta \theta}{\beta(1-\beta)}t}, \nonumber\\
 k(t)=k_0e^{-\frac{(\rho-\delta)(1-\beta)-\delta \theta}{\beta(1-\beta)}t},\nonumber\\
  u(t)=\bar{u}=\frac{\bigg(\rho-\delta(1-\beta
+ \theta)\bigg)(1-\beta)}{\delta \beta
(1-\beta+\theta)},\label{(sol1a)} \\
 h(t)=h_0 e^{-\frac{(\rho-\delta)(1-\beta)-\delta \theta}{\beta(1-\beta+\theta)}t}, \nonumber\\
\lambda(t)=c_0^{-\beta}e^{\frac{(\rho-\delta)(1-\beta)-\delta \theta}{(1-\beta)}t},\nonumber\\
 \mu(t)= \big(\frac{1-\beta+\theta}{1-\beta}\big)h_0^{\frac{\theta}{1-\beta}}c_1 e^{\frac{[(\rho-\delta)(1-\beta)-\delta \theta](\beta-\theta)}{\beta(1-\beta+\theta)}t},\nonumber\\
 z^*=\bigg(\frac {\delta +\pi}{\beta \gamma}
  \bigg)^{\frac{1}{1-\beta}} \nonumber,\end{eqnarray}
  provided $ \delta (1-\beta+\theta) <\rho,
c_0=\bigg(\frac{\gamma }{c_1 \delta (1-\beta+\theta)\bar{z}^{ \beta}
}\bigg)^{\frac{1}{\beta}},
 \frac{c_0}{k_0}=\frac{\rho+\pi(1-\beta)}{\beta}>0,h_0=\frac{\bar{z} k_0}{\bar{u}}$.
 The second closed-form solution for all variables is obtained by taking
$\sigma=\beta$ in {(\ref {(sol2)})} and is given by
 \begin{eqnarray} c(t)=c_0
z_0e^{-\frac{(\rho-\delta)(1-\beta)-\delta \theta}{\beta(1-\beta)}t}z(t)^{-1}, \nonumber\\
k(t)=k_0 z_0e^{-\frac{(\rho-\delta)(1-\beta)-\delta
\theta}{\beta(1-\beta)}t}z(t)^{-1}
,\nonumber\\
 h(t)=h_0  e^{-\frac{(\rho-\delta)(1-\beta)-\delta \theta}{\beta(1-\beta+\theta)}t}, \label{(sol2a)} \\
 u(t)=\bar{u}=\frac{\bigg(\rho-\delta(1-\beta
+ \theta)\bigg)(1-\beta)}{\delta \beta
(1-\beta+\theta)},\nonumber\\
 \lambda(t)=c_0^{-\beta}z_0^{-\beta}e^{\frac{(\rho-\delta)(1-\beta)-\delta \theta}{(1-\beta)}t}z(t)^{\beta},\nonumber\\
 \mu(t)= \big(\frac{1-\beta+\theta}{1-\beta}\big)h_0^{\frac{\theta}{1-\beta}}c_1 e^{\frac{[(\rho-\delta)(1-\beta)-\delta \theta](\beta-\theta)}{\beta(1-\beta+\theta)}t}
\nonumber,\end{eqnarray}

provided
 \begin{eqnarray}
\delta (1-\beta+\theta) <\rho, c_0=\bigg(\frac{\gamma }{c_1 \delta
(1-\beta+\theta)\bar{z}^{ \beta}
}\bigg)^{\frac{1}{\beta}},\nonumber\\
 \frac{c_0}{k_0}=\frac{\rho+\pi(1-\beta)}{\beta}>0,h_0=\frac{z_0 k_0}{\bar{u}}\nonumber\\
z(t)=
 \frac{\bar{z}z_0}{[(\bar{z}^{1-\beta}-z_0^{1-\beta})e^{-\frac{(1-\beta)(\delta
+\pi)+\delta \theta}{\beta}t}+z_0^{1-\beta}]^{\frac{1}{1-\beta}}
}\nonumber.\end{eqnarray}

It is important to mention here that the closed-form solution for
all the variables obtained by taking $\sigma=\beta$ in {(\ref
{(sol2a)})}  is the same as given in {(\ref {(sol3)})}.

In the previous literature, Ruiz-Tamarit \cite{ramon} established
one closed-form solution which was equal to our first solution
{(\ref {(sol1a)})} by utilizing the classical approach. Our
methodology has established two closed-form solutions ({(\ref
{(sol1a)})} and {(\ref {(sol2a)})} for the $\sigma=\beta$  case and
one of these solutions {(\ref {(sol2a)})} is completely new to the
literature.

\subsection{Convergence to balanced growth path}
For solution {(\ref {(sol1a)})}, the growth rates of per capita
consumption $c$, physical capital $k$, human capital $h$, the
fraction of labor allocated to the production of physical capital
$u$, costate variables $\mu$ and $\lambda$ take the following forms
(after some simplifications):
\begin{eqnarray}
 \frac{\dot c}{c}=\frac{(\delta -\rho)(1-\beta)+\delta \theta}{\beta(1-\beta)}, \nonumber \\
\frac{\dot k}{k}=\frac{(\delta -\rho)(1-\beta)+\delta
\theta}{\beta(1-\beta)},\nonumber \\
 \frac{\dot
{h}}{h}=\frac{(\delta -\rho)(1-\beta)+\delta
\theta}{\beta(1-\beta+\theta)},\nonumber \\
 \frac{\dot
u}{u}=0,\label{(gr1a)}\\
 \frac{\dot
{\lambda}}{\lambda}=\frac{(\rho-\delta)(1-\beta)-\delta
\theta}{1-\beta},\nonumber \\
 \frac{\dot
{\mu}}{\mu}=\frac{((\rho-\delta)(1-\beta)-\delta \theta
)(\beta-\theta)}{\beta(1-\beta+\theta)}. \nonumber
\end{eqnarray}

For solution {(\ref {(sol2a)})}, the growth rates of the per capita
consumption $c$, physical capital $k$, human capital $h$, the
fraction of labor allocated to the production of physical capital
$u$, costate variables $\mu$ and $\lambda$ take the following forms
(after some simplifications):
\begin{eqnarray}
 \frac{\dot c}{c}=\frac{(\delta -\rho)(1-\beta)+\delta \theta}{\beta(1-\beta)}-\frac{\dot z}{z}, \nonumber \\
\frac{\dot k}{k}=\frac{(\delta -\rho)(1-\beta)+\delta
\theta}{\beta(1-\beta)}-\frac{\dot z}{z},\nonumber \\
 \frac{\dot
{h}}{h}=\frac{(\delta -\rho)(1-\beta)+\delta
\theta}{\beta(1-\beta+\theta)},\nonumber \\
 \frac{\dot
u}{u}=0,\label{(gr1b)}\\
 \frac{\dot
{\lambda}}{\lambda}=\frac{(\rho-\delta)(1-\beta)-\delta
\theta}{1-\beta}+\beta\frac{\dot z}{z},\nonumber \\
 \frac{\dot
{\mu}}{\mu}=\frac{((\rho-\delta)(1-\beta)-\delta \theta
)(\beta-\theta)}{\beta(1-\beta+\theta)}, \nonumber
\end{eqnarray}
where $ \frac{\dot z}{z}$ is given in {(\ref {(grz)})}.

In both solutions we find the following: (i) the growth rates of
consumption $c$ and the physical capital $k$ decrease over time and
approach $\frac{(\delta -\rho)(1-\beta)+\delta
\theta}{\beta(1-\beta)}$ as $t\mapsto\infty$;  (ii) the growth rate
of human capital decreases over time and approaches $\frac{(\delta
-\rho)(1-\beta)+\delta \theta}{\beta(1-\beta+\theta)}$  as
$t\mapsto\infty$;  (iii) the growth rates of the per capita
consumption $c$, physical capital $k$ and human capital $h$ are
positive provided $(\delta -\rho)(1-\beta)+\delta \theta>0$; (iv)
the growth rate of the fraction of labor allocated to the production
of physical capital $u$ approaches zero as $t\mapsto\infty$; (v) the
growth rate of costate variable $\lambda$ converges to
$\frac{(\rho-\delta)(1-\beta)-\delta \theta}{1-\beta}$ as
$t\mapsto\infty$ and is negative; (vi) the growth rate of the
costate variable $\mu$ converges to
$\frac{((\rho-\delta)(1-\beta)-\delta \theta
)(\beta-\theta)}{\beta(1-\beta+\theta)}$ as $t\mapsto\infty$; and
(vii) the growth rate of the costate variable $\mu$ is negative when
$\beta>\theta$ and is positive when $\beta<\theta$.

Moreover, it is straight forward to check that both closed-form
solutions {(\ref {(sol1a)})} and {(\ref {(sol2a)})} satisfy all
properties of the BGP stated in Proposition 1.

\section{Conclusions}
The Lucas-Uzawa model with externalities is one of the fundamental
models of endogenous growth and is used to explain how economies
grow over time as they accumulate human capital. The solutions to
this model in the literature focus on how the key parameters in the
model, such as consumption, capital, and human capital, behave as
they progress on the balanced growth path towards equilibrium or
their long run equilibrium values.

But from a policymaker's perspective, the true benefit of such a
model would be to see how an economy behaves outside of equilibrium
in order to be able to design policies during an economy's
transition from the short-run (when it may be in disequilibrium) to
the long run. The issue that arises in much of the previous
literature is that most of methodologies used to derive solutions
outside of the equilibrium have yielded unique solutions. These
unique solutions are important but do not adequately explain how
economies starting from the same point can end up at distinctly
different points in their growth trajectories as they progress
towards a long run equilibrium.

We attempt to address this issue by using a newly developed
methodology to derive multiple solutions for the Lucas-Uzawa model
with externalities.  First, we began by transforming the Lucas-Uzawa
model with externalities into  the basic Lucas-Uzawa model. Then we
derived three closed-form solutions which are completely new to the
literature from the results in \cite{naz2017}. The closed-from
solutions presented in {(\ref {(sol1)})},  {(\ref {(sol2)})} and
{(\ref {(sol3)})} hold for fairly general values of the parameters.
 Also, while the previous literature, Hiraguchi
\cite{hir} derived only one solution for the model without parameter
restrictions which was similar to {(\ref {(sol3)})}, and in this
solution $F(t)$ and $G(t)$ were computed in terms of the
hypergeometric functions. Hiraguchi \cite{hir} also claimed that
there existed no other solutions for the Lucas-Uzawa model with
externalities.  Contrary to this claim and for the first time in the
literature, we have established multiple closed-form solutions for
the Lucas-Uzwa model with externalities in the case where there are
no parameter restrictions. Moreover, for $\sigma=\beta$, only one
solution {(\ref {(sol1)})} was established for the model
\cite{ramon}, we show that multiple solutions ({(\ref {(sol1a)})}
and {(\ref {(sol2a)})}) exist in this case.

Our methodology yields completely new multiple closed-form solutions
for the Lucas-Uzawa model with externalities and we use these
solutions to derive the growth rates for all of the variables in the
model. Our results are important on multiple levels: First, our
solutions establish multiplicity in terms of closed-form solutions
which has not been found in any of the previous economic growth
literature. This is important for our understanding of the
differences across countries in economic growth.  Second, one of our
solutions yields static growth rates whereas the other solutions
yield dynamic growth rates and in the long run all these growth
rates reach the same static value.  This means that these countries
that experience differences in growth over the short run can still
converge to the same long run growth rates.

\end{document}